\documentclass[12pt]{article}
\usepackage{times}
\usepackage{amsmath,amsxtra,amsthm,amssymb,makeidx,graphics}
\usepackage{latexsym,color}
\usepackage{tikz}
\usepackage{float}
\usepackage{euscript}

\theoremstyle{plain}
\newtheorem{theorem}{Theorem}[section]
\newtheorem{lemma}[theorem]{Lemma}
\newtheorem{corollary}[theorem]{Corollary}

\theoremstyle{definition}
\newtheorem{remark}[theorem]{Remark}
\newtheorem{definition}[theorem]{Definition}
\newtheorem{example}[theorem]{Example}
\newtheorem{exercise}[theorem]{Exercise}

\theoremstyle{remark}

\renewcommand{\emptyset}{\varnothing}
\newcommand{\ncom}{\newcommand}

\newcommand{\cB}{{\cal B}}

\newcommand{\cI}{{\cal I}}

\newcommand{\cL}{{\cal L}}
\newcommand{\cM}{{\cal M}}

\newcommand{\cR}{{\cal R}}

\newcommand{\cU}{{\cal U}}

\newcommand{\R}{{\mathbb R}}
\newcommand{\C}{{\mathbb C}}

\newcommand{\F}{{{\mathbb F}}}

\newcommand{\qb}[2]{{{ {{#1}\choose {#2}}_q }}}
\newcommand{\bin}[2]{{{ {#1}\choose {#2}}} }

\newcommand{\cspan}{{\mbox{span}}}

\newcommand{\rank}{{\mbox{rank}}}

\newcommand{\row}{{\cal{RE}}}

\newcommand{\ins}{{\mbox{ins}}}
\newcommand{\dele}{{\mbox{del}}}
\newcommand{\set}{{\mbox{set}}}
\newcommand{\subbset}{{\mbox{subset}}}
\newcommand{\sspace}{{\cal{RS}}}

\ncom{\n}[1]{{{ \|{#1}\|     }}}
\ncom{\ns}{\normalsize}
\ncom{\la}{\lambda}
\ncom{\bm}{\boldmath}
\ncom{\noi}{\noindent}
\ncom{\bq}{\begin{equation}}
\ncom{\eq}{\end{equation}}  
\ncom{\beqn}{\begin{eqnarray*}}
\ncom{\eeqn}{\end{eqnarray*}}  
\ncom{\ba}{\begin{array}}
\ncom{\ul}{\underline}   

\ncom{\ea}{\end{array}}
\ncom{\beq}{\begin{eqnarray}}
\ncom{\eeq}{\end{eqnarray}}  
\ncom{\nno}{\nonumber}
\ncom{\hs}{\mbox{\hspace{.25cm}}}
\ncom{\rar}{\rightarrow}
\ncom{\lrar}{\longrightarrow}
\ncom{\Rar}{\Rightarrow}
\ncom{\noin}{\noindent} 
\ncom{\bc}{\begin{center}}
\ncom{\ec}{\end{center}}  
\ncom{\sz}{\scriptsize}   
\ncom{\fpd}{\Phi(\pi^{'})}
\ncom{\fp}{\Phi(\pi) }
\ncom{\nk}{\left< \begin{array}{c}
                       n\\k \end{array} \right>}
\ncom{\nd}{1^{'},2^{'},\cdots,n^{'}}
\ncom{\de}{\bigtriangleup (F_{2n},\leq)}
\ncom{\del}{\bigtriangleup}
\ncom{\cov}{<\!\!\!\!\cdot }
\ncom{\bt}{\begin{theorem}}
\ncom{\bcon}{\begin{con}}
\ncom{\et}{\end{theorem}}
\ncom{\econ}{\end{con}}
\ncom{\bl}{\begin{lemma}}
\ncom{\el}{\end{lemma}}  
\ncom{\bco}{\begin{corollary}} 
\ncom{\ds}{\displaystyle}
\ncom{\eco}{\end{corollary}}   
\ncom{\bp}{\begin{pro}}  
\ncom{\ep}{\end{pro}}    
\ncom{\bex}{\begin{example}}
\ncom{\eex}{\end{example}}  
\ncom{\bexr}{\begin{exercise}}
\ncom{\eexr}{\end{exercise}}  
\ncom{\bprob}{\begin{problems}}
\ncom{\eprob}{\end{problems}}  
\ncom{\bd}{\begin{definition}}
\ncom{\ed}{\end{definition}}  
\ncom{\brm}{\begin{remark}}   
\ncom{\erm}{\end{remark}}     
\ncom{\bal}{\begin{Algorithm}}
\ncom{\eal}{\end{Algorithm}}  
\ncom{\ol}{\overline}
\ncom{\wh}{\widehat} 

\ncom{\pf}{\noi {\bf Proof.  }}
\ncom{\eprf}{\noi {$\Box$}}
\ncom{\be}{\begin{enumerate}} 
\ncom{\ee}{\end{enumerate}}   
\ncom{\seq}{\subseteq}

\ncom{\zr}{\bf\textcolor{red}}
\ncom{\zb}{\bf\textcolor{blue}}
\ncom{\zg}{\bf\textcolor{green}}
\ncom{\zm}{\bf\textcolor{magenta}}

\definecolor{gold}{rgb}{0.85,.66,0}
\definecolor{gb}{rgb}{0, .5,.5}
\definecolor{rb}{rgb}{0.5, 0,.5}
\definecolor{Pink}{rgb}{1,0.75, 0.8}
\ncom{\zgb}{\bf\textcolor{gb}}
\ncom{\zrb}{\bf\textcolor{rb}}

\makeatletter
\newcommand{\single}{\let\CS=\@currsize\renewcommand{\baselinestretch}{1.5}\tiny\CS}
\newcommand{\oneandahalfspacing}{\let\CS=\@currsize\renewcommand{\baselinestretch}{1.5}\tiny\CS}
\newcommand{\doublespacing}{\let\CS=\@currsize\renewcommand{\baselinestretch}{1.6}\tiny\CS}
\newcommand{\double}{\let\CS=\@currsize\renewcommand{\baselinestretch}{3}\tiny\CS}

\def\notarro{{{\hbox{{\hspace*{-.02in}}$\rightarrow$ } }
{\hbox{$\!\!\!\!\!\!\!$}}{\raise 0.2ex 
\hbox{$\scriptscriptstyle{/}$}}}\hspace{.06in}}
\def \*{^{\mbox{$*$}}}

\newcommand\bnota{\begin{nota} }                         
\newcommand\enota{\end{nota} }                           

\newcommand{\beano}{\begin{eqnarray*}}
\newcommand{\eeano}{\end{eqnarray*}}

\begin{document}

\title{\bf{\Large{Subspaces, subsets, and Motzkin paths}}}
 
\author{  \normalsize{ \bf{Jonathan D. Farley}} \\
            \em{Department of Mathematics} \\
            \em{1700 E. Cold Spring Lane} \\
            \em{Morgan State University}\\
            \em{Baltimore, MD 21251, USA}\\
            email: lattice.theory@gmail.com  
\and        \normalsize{\bf{Murali K. Srinivasan}}\\
             \em{Department of Mathematics} \\
             \em{Indian Institute of Technology, Bombay} \\
             \em{Powai, Mumbai 400076, INDIA} \\
             email: murali.k.srinivasan@gmail.com   }
\date{}
\maketitle
\begin{center}
{\small{\em To the memory of my father in law Somasundaram Lakshmipathy (MKS)}}
\end{center}

\begin{abstract}

Let $\cM(n)$ denote the set of all Motzkin paths from $(0,0)$ to $(n,0)$.
For each $P\in \cM(n)$ we define a statistic $w(P,q)$, the weight of $P$.
Let $|P|$ denote the number of down steps in $P\in\cM(n)$. 
Let $B_q(n)$ denote the projective geometry
(= poset of subspaces of an $n$-dimensional vector space over
$\F_q$).

We define a map from $B_q(n)$ to $\cM(n)$ and show that, for $P\in \cM(n)$,
the inverse image of $P$ consists of a disjoint union of 
$(q-1)^{|P|} w(P,q)$ symmetric Boolean subsets
in $B_q(n)$, all with minimum rank $|P|$ and maximum rank $n-|P|$. This yields
an explicit symmetric Boolean decomposition of the projective geometry
and gives a poset theoretic interpretation to the identity
\beqn
 \qb{n}{k} &=& \sum_{P \in \cM(n)}
                       (q-1)^{|P|} w(P,q)
                   \bin{n-2|P|}{k-|P|}.
\eeqn
\end{abstract}
     
\noi

\section{Introduction}

This paper combines two results:

(i) The solution, by Vogt and Voigt {\bf\cite{vv}}, of 
Greene and Kleitman's {\bf\cite{gk}} problem of constructing an 
explicit symmetric chain decomposition of the subspace lattice.  

(ii) The symmetric expansion, from the 
paper {\bf\cite{ds}}, 
of the $q$-binomial coefficient in terms
of the binomial coefficients with summands indexed by involutions, 
obtained by iterating the Goldman-Rota recurrence {\bf\cite{gr}} for the  
number of subspaces of a finite vector space.

These two results are closely related. The idea is to consider both
of them through the lens of Motzkin paths.
Biane {\bf\cite{b1}} gave a map from involutions to Motzkin paths
together
with a simple formula for the cardinality of the inverse image of 
a Motzkin path (we learnt about this map from {\bf\cite{bbs}}).
We show that this formula has a $q$-analog and using this we 
rewrite the expansion from {\bf\cite{ds}}, with the summands now indexed
by Motzkin paths. Then we give the resulting 
identity a poset theoretic interpretation by defining a map from 
subspaces to Motzkin paths and showing that 
\begin{itemize}
\item the cardinality 
of the inverse image of every
Motzkin path agrees with the formula and that 
\item the inverse images are 
a disjoint union of symmetric Boolean subsets. 
\end{itemize}
The definition
of this map is remarkably simple and is inspired by Biane's map.
The introduction of Motzkin paths and the map from subspaces to Motzkin paths
reveals the underlying structure in the subspace 
lattice found by {\bf\cite{vv}}. 

Let us now state our results precisely.

Let $P$ be a finite graded poset of rank $n$   
with rank function $r$.
For
$0 \leq k \leq n$, let $N_k$ denote the number of elements of $P$ of rank
$k$.
We say that the sequence of elements $(x_1, x_2, \ldots ,x_h)$ of $P$  form a {\it
symmetric
chain} if $x_{i+1}$ covers $x_{i}$ for every $i < h$ and $r(x_1) +
r(x_h)=n$ if $h\geq 2$ or $2r(x_1)=n$ if $h=1$. 
A {\it symmetric chain decomposition} (SCD) of $P$ is a covering of $P$ by    
pairwise disjoint symmetric chains.  Let $B(n)$ denote the {\em Boolean algebra}, i.e.,
the graded
poset, under inclusion, of all
subsets of $[n]=\{1,2,\ldots ,n\}$, 
where the rank of a subset is its cardinality.
We say that a subset $Q \subseteq P$ is {\it 
symmetric Boolean} if 

\noi
(i) $Q$, under the induced order, has a minimum element, say $z$, and a
maximum
element, say $z^{\prime}$. \\
(ii) $Q$ is order isomorphic to $B(r(z^{\prime})-r(z))$. \\
(iii) $r(z^{\prime}) + r(z) =n$. 

A {\it symmetric Boolean decomposition} (SBD) of $P$ is a covering of $P$ by
pairwise disjoint symmetric Boolean subsets.
De Bruijn, Tenbergen, and Kruyswijk {\bf\cite{btk}} inductively constructed a symmetric 
chain decomposition of $B(l)$, 
for $l \geq 0$ (and, more generally, for chain products). Several authors have given an explicit
version of this SCD, see {\bf\cite{gk}}. It follows that if $P$ admits
a SBD, then it has a SCD. Moreover, an explicit construction of a SBD immediately
yields an explicit construction of a SCD.

The existence of a SBD gives a symmetric expansion 
of the rank numbers of $P$ in terms
of the binomial coefficients. Let $P = Q_1 \cup Q_2 \cup \cdots
\cup Q_t $  
(disjoint union) be a SBD of $P$. Let $z_i$ (respectively $z_i^{\prime})$
denote the minimum (respectively, maximum) element of $Q_i$, $i=1,2, \ldots ,t$.
Since $Q_i$ is order isomorphic to $B(r(z_i^{\prime})-r(z_i))$ and
$r(z_i^{\prime}) + r(z_i) = n$ we have
\beqn
\label{rid1}
N_k =  \sum_{i=1}^{t} {r(z_i')-r(z_i) \choose k-r(z_i)} = 
\sum_{i=1}^{t} {n-2r(z_i)\choose k-r(z_i)},
\eeqn
and, summing over $k$, we get
\beqn
|P|  =  \sum_{i=1}^{t} 2^{n-2r(z_i)}.
\eeqn
Let $q$ be a prime power and let $B_q(n)$ denote the {\em (finite) projective geometry}
, i.e., the graded poset of subspaces,
under inclusion, of the $n$-dimensional vector space $\F^n_q$ over $\F_q$ 
(we think of $\F^n_q$ as row vectors). The rank of a subspace is its dimension.
The number of elements of rank $k$ in
$B_q(n)$
is the $q$-binomial coefficient $\qb{n}{k}$. 
Griggs {\bf\cite{g}} proved, using network flow techniques,
that
$B_q(n)$ has a SCD and Greene and Kleitman {\bf\cite{gk}} 
asked for an explicit construction.
Bj\"{o}rner asked (see Exercise 7.36 in {\bf\cite{b2}}) 
whether $B_q(n)$ has a SBD (the exercise actually concerns the related concept
of Boolean packings (defined by relaxing condition (iii) in the definition of a SBD
above to $r(z^{\prime}) + r(z) \geq n$) 
for all finite geometric lattices but for the subspace lattice this is
just SBD. 
Boolean packings are very easy to construct for set partition lattices,
see {\bf\cite{s1}}). An unsuccessful attempt, in {\bf\cite{ds}}, to construct
an explicit SBD of $B_q(n)$ led to a symmetric expansion of $\qb{n}{k}$ of the form
above, with the summands indexed by involutions (see Theorem \ref{ds} in Section 2).

Let $\cM(n)$ denote the set of Motzkin paths from $(0,0)$ to $(n,0)$ , i.e., 
lattice paths from $(0,0)$ to $(n,0)$ with steps $(1,0)$ ({\em horizontal 
step}), $(1,1)$ ({\em up step}), and $(1,-1)$ ({\em down step}), never going 
(strictly) below the $x$-axis. Define the {\em height}
of a down step $(i,j+1)$ to $(i+1, j)$ to be $j+1$. 
Define the {\em weight} of 
\begin{itemize}
\item an up  step to be $1$.
 
\item the weight of a horizontal step $(i,j)$ to $(i+1,j)$ to be $q^j$.

\item the weight of a down  step $(i,j+1)$ to $(i+1,j)$ 
to be $q^j + q^{j+1} + \cdots + q^{2j}$. So the weight of a down step 
is the height when $q=1$. 
\end{itemize}
The {\em weight} $w(P,q)$ of a Motzkin path $P$ is the product of the 
weights of the steps of $P$ and let $|P|$ denote the number of down 
steps of $P$. We shall write a Motzkin path in $\cM(n)$ as $s_1s_2\cdots s_n$,
where each $s_i\in\{U,D,H\}$. For instance, the Motzkin path
$UHDHUUHDD \in \cM(9)$ pictured below

\begin{figure}[h]
    \centering
    \begin{tikzpicture}[scale=0.5]
        \draw[step=1cm,gray,very thin] (0,0) grid (9,4);
        \draw[->,thick] (0,0) -- (9.5,0) node[anchor=west] {};
        \draw[->,thick] (0,0) -- (0,4.5) node[anchor=south] {};
        \draw[ultra thick] (0,0) -- (1,1) -- (2,1) -- (3,0) -- (4,0) -- (5,1) -- (6,2) -- (7,2) -- (8,1) -- (9,0);
        \node[below] at (0,0) {\tiny (0,0)};
        \node[below] at (9,0) {\tiny (9,0)};
    \end{tikzpicture}
    \caption{\( UHDHUUHDD \)}
    \label{fig:motzkin_path}
\end{figure}
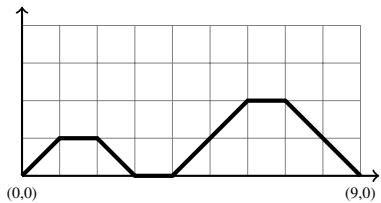

has weight $1\cdot q \cdot 1 \cdot 1 \cdot 1 \cdot 1 \cdot q^2\cdot (q+q^2)
\cdot 1 = q^3(q+q^2)$.

\bt \label{fs}
We have the following $q$-binomial expansion
\beq \label{qbi}
 \qb{n}{k} &=& \sum_{P \in \cM(n)} 
                       (q-1)^{|P|} w(P,q)
                   {\bin{n-2|P|}{k-|P|}}.
\eeq
\et

In Section 2 we give a manipulatorial proof of Theorem \ref{fs}. We first 
recall, in Theorem \ref{ds},
the identity from {\bf\cite{ds}} (similar
to (\ref{qbi}) above except that the sum is over involutions) that follows 
from iterating the Goldman-Rota recurrence. Then we recall Biane's map 
from involutions to Motzkin paths, together with a simple formula for the
cardinality of an inverse image. We show that this formula has a $q$-analog  
and collect terms in Theorem \ref{ds}
in accordance with this $q$-analog, proving
Theorem \ref{fs}.

In Section 3 we define a map $\Psi:B_q(n)\rar \cM(n)$, a vector space
analog of Biane's map, and 
prove the following result which gives a poset theoretic interpretation to
Theorem \ref{fs}.

\bt \label{main} For $P\in \cM(n)$, $\Psi^{-1}(P)$ is a disjoint
union of $(q-1)^{|P|}w(P,q)$ symmetric Boolean subsets, all with minimum rank
$|P|$ and maximum rank $n-|P|$.
\et

There is an alternative way of viewing 
the map $\Psi$. 
The usual Gauss elimination method classifies the columns 
of a matrix in row reduced echelon form as pivotal/nonpivotal. 
Motivated by the approach
in {\bf\cite{vv}} we introduce a different classification into
essential/inessential columns. Combining these two gives 
a four fold classification 
of the columns: each column can be pivotal/nonpivotal 
and essential/inessential. 
Theorem \ref{cqbiane} characterizes the map $\Psi$ in terms of 
this classification and this is then used in the proof of Theorem \ref{main}.
Using Theorem \ref{main} we also give a simpler and more elegant description
of the SCD of $B_q(n)$ constructed in {\bf\cite{vv}}. 

Finally, in Section 4, we state a question concerning a possible 
relation between the set theoretic and linear analogs of  
symmetric chains in $B(n)$ and $B_q(n)$.

{\bf{Remark}} In the paper {\bf\cite{ch}} by Coopman and Hamaker
Theorem \ref{fs} and Theorem \ref{ds} (from Section 2)
are combined into a single equation numbered (4.2) and stated without proof. The authors offer 
the following comment: ``The proof of Eq. (4.2) follows from
generating function manipulations, and it would be interesting to give a direct
combinatorial proof using Motzkin paths". Theorem \ref{main} presents
such a direct combinatorial proof. 
We have given a detailed proof of Theorem \ref{fs} above since exactly the same
pattern of argument occurs in part of the proof of Theorem \ref{main} 
(see the similarities in the proof of Lemma \ref{qa} and Theorem \ref{bpb}).

The title of the present paper was inspired by {\bf\cite{k}} which used row reduced echelon forms to give a direct connection between subspaces and partitions.

\section{Involutions and Motzkin Paths}

Let $G_q(n)=\sum_{k=0}^n \qb{n}{k}$, the {\em Galois numbers},
denote the number of subspaces of $\F^n_q$. The starting point of this paper
is the Goldman-Rota recurrence {\bf\cite{gr,kc}}
\beq \label{gr}
G_q(n+1) &=& 2G_q(n) + (q^n-1)G_q(n-1),\;\; G_q(0)=1,G_q(1)=2.
\eeq
Unfolding the recurrence we can expand $G_q(n)$ in powers of two. For instance,
$$G_q(2)=2^2 +(q-1)2^0,\;\;G_q(3)=2^3 +(q-1+q^2-1)2^1.$$
To get the coefficients in the 
expansion we write $q^n-1 =(q-1)(1+q+q^2+\cdots +q^{n-1})$
and rewrite the recurrence as
\beqn
\lefteqn{G_q(n+1) = 2G_q(n) +}\\
&& (q-1)G_q(n-1) +q(q-1)G_q(n-1) + \cdots +q^{n-1}(q-1)G_q(n-1)
\eeqn
The right hand side has one ocurrence of $G_q(n)$ (with a coefficient of 2)
and $n$ occurrences
of $G_q(n-1)$ (with coefficients $q^i(q-1),\;i=0,1,\ldots ,n-1$). 
Informally speaking, this is the same recursive structure
as that of involutions in the symmetric group $S(n+1)$ (either $n+1$ is a fixed
point or it is paired with the letter $i$, where  $i\in\{1,2,\ldots ,n\}$).
This suggests that there is an expansion of $G_q(n)$ in terms of 
powers of 2 with the summands
indexed by involutions in $S(n)$ and the coefficients 
given by some 
statistic on involutions. We now recall this identity.

Let $\cI(n)$ denote the set of all involutions in $S(n)$.
We write a 2-cycle in $S(n)$ as $[i,j]$, with $i < j$. For a 2-cycle
$[i,j]$, we call $i$ the {\it initial point} and $j$ the {\it terminal
point}.  The {\it span} of $[i,j]$ is defined as
$\cspan([i,j])=j-i-1$. A pair $\{[i,j], [k,l]\}$ of disjoint
2-cycles is said to be a {\it crossing} if $i < k < j < l$ or $k < i < l < j$
(see Figure 2).
We write involutions in their
{\it standard form} by listing the 2-cycles 
in increasing order of initial points.

\begin{figure}[h]
\includegraphics[scale=.8]{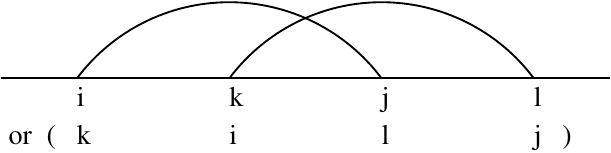}
\caption{A Crossing}
\end{figure}

Let $\delta$ be an involution in $\cI(n)$. The number of 2-cycles 
in $\delta$ is denoted  
by $|\delta|$. The {\it crossing number} of $\delta$, denoted $c(\delta)$,
is the number of pairs of 2-cycles of $\delta $ that are crossings.
Define the {\it weight} of $\delta$, denoted by $w(\delta)$, as follows
$$
w(\delta) = \left(\displaystyle\sum_{[i,j]}\cspan([i,j])
                                   \right) - c(\delta),
$$
where the sum is over all 2-cycles in $\delta$.
\bex
Let $\delta = [1,8][2,6][3,9][4,7] \in \cI(9)$. Represent
$\delta$ as shown in Figure 3. Observe that there are 3
crossings. Thus $w(\delta) =
(8-1-1)+(6-2-1)+(9-3-1)+(7-4-1)-3=13$.
\begin{figure}[h]
\includegraphics[scale=.8]{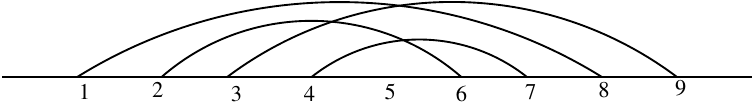}
\caption{Crossing Number}
\end{figure}

\eex

The following result was proved in {\bf\cite{ds}}.
\bt \label{ds}
For $0\leq k\leq n$ we have
\beq \label{mkst}
\qb{n}{k} &=&  \sum_{\delta \in
                   \cI(n)} (q-1)^{|\delta|}q^{w(\delta)}
                   \bin{n-2|\delta|}{k-|\delta|}.
\eeq  
Summing over $k$ we obtain
\beqn G_q(n) &=&  \sum_{\delta \in
                   \cI(n)} (q-1)^{|\delta|}q^{w(\delta)}
                    2^{n-2|\delta|}.
\eeqn  
\et
For example, $\qb{5}{k}$  equals
$$ 
\bin{5}{k} + (q-1)(4+3q+2q^{2}+q^{3})\bin{3}{k-1}
+ (q-1)^{2}(3+4q+4q^{2}+3q^{3}+q^{4})\bin{1}{k-2}.
$$

We do not have a poset theoretic interpretation of Theorem \ref{ds} in the manner
of Theorem \ref{main}. Instead, we shall rewrite the 
identity (\ref{mkst}) in terms of Motzkin paths.

Define a map ({\bf\cite{b1,bbs}})
$$\cB: \cI(n)\rar \cM(n)$$
as follows:
write $\delta\in \cI(n)$ in standard form as $\delta=[i_1,j_1][i_2, j_2]
\cdots [i_k,j_k]$. Then  $\cB(\delta)=s_1s_2 \cdots s_n$, where
\beqn
s_t &=& \left\{ \ba{ll}
                   U,&\;\;t\in \{i_1,i_2,\ldots ,i_k\},\\
                   D,&\;\;t\in \{j_1,j_2,\ldots ,j_k\},\\
                   H,&\;\;i\in [n] - \{i_1,\ldots ,i_k,j_1,\ldots ,j_k\}.\ea \right.
\eeqn
Note that
\beq \label{invmp}
|\cB(\delta)| = |\delta|,\;\;\delta\in \cI(n).
\eeq
\bex The involution $[1,6][3,5]\in \cI(6)$ is represented by
the Motzkin path (Figure 4)

\begin{figure}[h]
    \centering
    \begin{tikzpicture}[scale=0.5]
        \draw[step=1cm,gray,very thin] (0,0) grid (7,4);
        \draw[->,thick] (0,0) -- (7.5,0) node[anchor=west] {};
        \draw[->,thick] (0,0) -- (0,4.5) node[anchor=south] {};
        \draw[ultra thick] (0,0) -- (1,1) -- (2,1) -- (3,2) -- (4,2) -- (5,1) -- (6,0);
        \node[below] at (0,0) {\tiny (0,0)};
        \node[below] at (6,0) {\tiny (6,0)};
    \end{tikzpicture}
    \caption{\(UHUHDD\)}
    \label{fig:motzkin_path}
\end{figure}

\eex

The following is clear (see {\bf\cite{b1,bbs}}).
\bl Let $P\in \cM(n)$. Then
$$|\cB^{-1}(P)|=\prod_{s}h(s),$$
where the product is over all down steps $s$ in $P$.
\el
The result above has a $q$-analog.
\bl \label{qa}
Let $P\in \cM(n)$. Then
$$\sum_\delta q^{w(\delta)}=w(P,q),$$
where the sum is over all $\delta\in \cI(n)$ with $\cB(\delta)=P$.
\el
\pf By induction on $n$, the case $n=1$ being clear. Let $P=s_1s_2\cdots s_n\in \cM(n)$.
Let $t\geq 0$ be the largest integer with 
$s_1=s_2=\cdots =s_t=U$. So $s_{t+1}= H\mbox{ or }D$.

Let $\delta\in \cI(n)$ satisfy $\cB(\delta)=P$. From $s_1=\cdots =s_t=U$ we deduce
that the standard form of $\delta$ begins
$[1,\cdot][2,\cdot]\cdots [t,\cdot]\cdots $. Now consider two cases:

(i) $s_{t+1}=H$: In this case $t+1$ will not appear in the standard form 
of $\delta$. Taking the standard form of $\delta$ and subtracting
1 from every integer $>t+1$ we get the standard form of an involution
$\delta'\in\cI(n-1)$. It is easy to see that
\beq \label{wi}
w(\delta)&=&t+w(\delta'),
\eeq 
since the point $t+1$ contributes 1 to the span of the $t$ arcs
in $\delta$ with initial points $1,\ldots ,t$ and $c(\delta)=c(\delta')$.

Let $P'$ be the Motzkin path $s_1\cdots s_ts_{t+2}\cdots s_n \in \cM(n-1)$. 
By the induction hypothesis, $q$-counting with respect to weight,
the involutions $\delta''\in\cI(n-1)$ with $\cB(\delta'')=P'$ gives 
$w(P',q)$. 
Taking the standard
form of any such involution and increasing by 1 all integers
$\geq t+1$ gives the standard form of an involution in $\cI(n)$ 
mapping to $P$ under the map $\cB$. Since $w(P,q)=q^tw(P',q)$
the assertion now follows from (\ref{wi}).

(ii) $s_{t+1}=D$: In this case the standard form of $\delta$
begins 
$$[1,\cdot][2,\cdot]\cdots [j-1,\cdot][j,t+1][j+1,\cdot]\cdots [t,\cdot]\cdots$$
for some $j\in \{1,\cdots ,t\}$.

Take the standard form of $\delta$, remove the 2-cycle $[j,t+1]$, subtract
1 from each of $j+1,\ldots ,t$, and subtract 2 from each of $t+2,\ldots ,n$ to
get the standard form of an involution $\delta'\in\cI(n-2)$. We claim that
\beq \label{wi1}
w(\delta)&=&t+j-2+w(\delta').
\eeq
This can be seen as follows. In $\delta$, the 2-cycle $[j, t+1]$ has span
$t-j$ and  also participates
in $t-j$ crossings (with the $t-j$ 2-cycles whose initial points are $j+1,\ldots ,t$). So these two contributions to $w(\delta)$ cancel. The number of
crossings between the other arcs of $\delta$ is the same as $c(\delta')$. 
When going from $\delta$ to $\delta'$, the span of every arc with initial point
$1,\ldots ,j-1$ decreases by 2, while the span of every arc with initial point
$j+1,\ldots ,t$ decreases by 1. So $w(\delta)=2(j-1)+t-j+w(\delta')$ proving the claim.

Let $P'$ be the Motzkin path $s_1\cdots s_{t-1}s_{t+2}\cdots s_n \in \cM(n-2)$.
We have 
\beq \label{wi2}
w(P,q)=(q^{t-1}+q^t+\cdots +q^{2t-2})w(P',q).
\eeq
By the induction hypothesis, $q$-counting with respect to weight,
the involutions $\delta''\in\cI(n-2)$ with $\cB(\delta'')=P'$ gives 
$w(P',q)$. 
Taking the standard
form of any such involution, 
adding 2 to each of $t,\ldots , n-2$ and, for some $j=1,\ldots ,t$,
adding 1 to each of $j,\ldots ,t-1$ and adding the 2-cycle $[j,t+1]$,
we get an involution in $\cI(n)$ mapping to $P$ under $\cB$. 
The assertion follows from (\ref{wi1}) and (\ref{wi2}). \eprf

{\bf{Proof of Theorem \ref{fs}}} Follows from Theorem \ref{ds}, (\ref{invmp}), and Lemma
\ref{qa}. \eprf

\section{Subspaces and Motzkin Paths}

We define the map $\Psi:B_q(n)\rar \cM(n)$ from the introduction and prove
Theorem \ref{main}.

Let $e_1,e_2,\ldots ,e_n$ denote the standard basis of $\F_q^n$ (row vectors).
Define a map
$$\cL:B_q(n)\rar B(n)$$
by 
\beqn
\cL(X) &=& \left\{j\in [n]\;:\; e_j + \sum_{i>j}\alpha_{ij}e_i\in X 
\mbox{ for some } \alpha_{ij}\in \F_q\right\}.
\eeqn
It is easy to see that $\cL$ is rank and order preserving. We can compute $\cL$
using Gauss elimination. A $k\times n$ matrix over $\F_q$ is said to be in 
{\em row reduced echelon form} (rref) (also called {\em Schubert normal form})
provided
\begin{itemize}
\item Every row is nonzero and the first nonzero entry (from the left)
in every row is 1. Let the first nonzero entry in row $i$ occur in column $p_i$.

\item $p_1<p_2<\cdots <p_k$.

\item Columns $p_1,p_2,\ldots ,p_k$ form the $k\times k$ identity matrix.
\end{itemize}
We call $p_1,\ldots ,p_k$ the {\em left pivotal columns}. It is well known
that every $k$-dimensional subspace $X$ of $\F^n_q$ is the 
row space of a unique
$k\times n$ matrix in rref which can be computed by Gauss elimination, 
more precisely left to right Gauss elimination, starting from
any matrix with row space $X$. We denote the unique $k\times n$ rref with
row space $X$ by $\row(X)$. Clearly, $\cL(X)=$ the left pivotal
columns of $\row(X)$.                     

Similarly, define a map
$$\cR:B_q(n)\rar B(n)$$
by 
\beqn
\cR(X) &=& \left\{j\in [n]\;:\; e_j + \sum_{i<j}\alpha_{ij}e_i\in X 
\mbox{ for some } \alpha_{ij}\in \F_q\right\}.
\eeqn
We can compute $\cR$ by a right to left variant of
Gauss elimination. 
We call elements of $\cR(X)$ the {\em right pivotal
columns} of $X$. The term {\em pivotal column} will mean left pivotal column.                     

We can now define the map $\Psi$. 
Note the similarity with Biane's map.
Let $X\in B_q(n)$. Define $\Psi(X)=s_1s_2 \cdots s_n$, where (below $\Delta$
denotes symmetric difference)
\beqn
s_i &=& \left\{ \ba{ll}
                   U,&\;\;i\in \cL(X)\setminus \cR(X),\\
                   D,&\;\;i\in \cR(X)\setminus \cL(X),\\
                   H,&\;\;i\in [n]\setminus (\cL(X) \triangle \cR(X)).\ea \right.
\eeqn
Informally speaking, one feels that $\cR(X)$ should lie to the right
of $\cL(X)$. 
This basic property of Gauss elimination is made precise in the following result.
\bt \label{qbiane}
For all $X\in B_q(n)$, we have $\Psi(X)\in \cM(n)$.
\et
Before proving Theorem \ref{qbiane} let us see an example. First we introduce 
some notation. For a matrix $M$ in rref we let $\sspace(M)$ be the row space of
$M$ and we set $\cL(M)=\cL(X)$, $\cR(M)=\cR(X)$, and $\Psi(M)=\Psi(X)$, where
$X=\sspace(M)$.
\bex \label{exmap}
Consider the following matrix $M$ in rref
$$\left[ \ba{cccccc} 1&a&0&b&0&0\\0&0&1&c&0&d\\0&0&0&0&1&e \ea \right]$$
where $a, c, e\in \F_q$ and $b, d\in \F_q^*$. Then $\cL(M)=\{1,3,5\}$,
$\cR(M)=\{4,5,6\}$ and $\Psi(M)=UHUDHD$ (Figure 5 below).

\begin{figure}[h]
    \centering
    \begin{tikzpicture}[scale=0.5]
        \draw[step=1cm,gray,very thin] (0,0) grid (7,4);
        \draw[->,thick] (0,0) -- (7.5,0) node[anchor=west] {};
        \draw[->,thick] (0,0) -- (0,4.5) node[anchor=south] {};
        \draw[ultra thick] (0,0) -- (1,1) -- (2,1) -- (3,2) -- (4,1) -- (5,1) -- (6,0);
        \node[below] at (0,0) {\tiny (0,0)};
        \node[below] at (6,0) {\tiny (6,0)};
    \end{tikzpicture}
    \caption{ \( UHUDHD \)}
    \label{fig:motzkin_path}
\end{figure}

\eex
To prove Theorem \ref{qbiane} it is useful to have a 
left to right characterization 
of the right pivotal columns $\cR(X)$. Since there is a 
bijection between $k$-dimensional subspaces and $k\times n$
rref's it should be possible to calculate $\cR(X)$ from the rref 
representing $X$
without having to run the right to left variant of Gauss elimination. 
We discuss this next.

Let $M$ be a $k\times n$ matrix in rref with columns $C[1], \ldots ,C[n]$
and rows $R[1], \ldots ,R[k]$. For $0\leq m \leq k$, let $C_m[j]$ denote
the column vector formed by the {\em{first}} $m$ components (from the top) 
of $C[j]$. For $0\leq r\leq n$,
let $R_r[i]$ denote the row vector formed 
by the {\em{last}} $r$ components (from the left) 
of
$R[i]$. Denote the pivotal columns of $M$ by $\{p_1 < p_2 < \cdots < p_k\}$.

Let $j\in \{1,\ldots ,n\}$. The {\em{section $S_j$ at $j$}} is the submatrix of $M$
defined as follows:

(i) $j$ is nonpivotal: Let $m$ be the unique integer with $p_m<j<p_{m+1}$. Then $S_j$ is the submatrix
of $M$ formed by the first $m$ rows and last $n-j$ columns of $M$, i.e., the rows
of $S_j$ are $R_{n-j}[1], R_{n-j}[2], \ldots ,R_{n-j}[m]$ and the columns of
$S_j$ are $C_m[j+1],\ldots ,C_m[n]$. We can picture $M$ as follows

$$\left[\ba{ccc} A &C_m[j] & N\\ 
         0 & 0 & B \ea \right]$$

where $A$ is $m\times (j-1)$, $N$ is $m\times (n-j)$, $B$ is $(k-m)\times (n-j)$, and $S_j=N$.

(ii) $j$ is pivotal: Let $j=p_m$. Then $S_j$ is the submatrix
of $M$ formed by the first $m$ rows and last $n-j$ columns of $M$, i.e., the rows
of $S_j$ are $R_{n-j}[1], R_{n-j}[2], \ldots ,R_{n-j}[m]$ and the columns of
$S_j$ are $C_m[j+1],C_m[j+2],\ldots ,C_m[n]$. We can picture $M$ as follows

$$\left[ \ba{ccc} A&0&N\\
                  0&1&R_{n-j}[m]\\
                  0&0&B \ea\right]$$

where $A$ is $(m-1)\times (j-1)$, $N$ is $(m-1)\times (n-j)$, $B$ is $(k-m)\times (n-j)$, and $S_j=\left[\ba{c}N\\R_{n-j}[m]\ea\right]$.

Note that $S_n$ is the empty matrix. We also define $S_0$ to be the empty matrix.

Column $j$ of $M$ is said to be {\em{essential}} if the following holds. We
consider two cases.

(i) $j$ is nonpivotal: Let $m$ be the unique integer with  
$p_m < j < p_{m+1}$. Then 
$$C_m[j] \not\in 
\cspan\{C_m[j+1], C_m[j+2], \ldots ,C_m[n]\},$$
i.e., $C_m[j]$ is not in the column space of $S_j$.

(ii) $j$ is pivotal: Let  $j=p_m$. Then 
$$R_{n-j}[m] \not\in 
\cspan\{R_{n-j}[1], R_{n-j}[2], \ldots ,R_{n-j}[m-1]\},$$
i.e., the last row of $S_j$ does not linearly depend on the other rows.

A column that is not essential is said to be {\em{inessential}}. So we have
four types of columns: pivotal and nonpivotal, essential and inessential.
\bex (i) Consider column 1. If it is nonpivotal then it is
inessential (since $S_1$ is the empty matrix). If column 1 is pivotal then it is inessential if and only if the first row is $e_1$.

(ii) Consider column $n$. It it is pivotal then it is inessential (since $S_n$ is the empty matrix). If column $n$ is nonpivotal then it is inessential if and only if it is the zero column.

(iii) In Example \ref{exmap} the inessential columns are $2,5$ of which
column $2$ is nonpivotal and column $5$ is pivotal. The essential columns are
$1,3,4,6$ of which $1,3$ are pivotal and $4,6$ are nonpivotal. 
\eex

The next result relates the definitions above to the map $\Psi$.
\bt \label{cqbiane}
Let $X\in B_q(n)$ and let $M=\row(X)$ be $k\times n$.
Let $\Psi(X)=s_1s_2\cdots s_n$. Then

(i) $s_j=H$ if and only if $j$ is inessential.

(ii) $s_j=U$ if and only if $j$ is essential and pivotal.

(iii) $s_j=D$ if and only if $j$ is essential and nonpivotal.
\et
\pf 
We shall use the notation for the rows, columns, and pivotal columns of $M$
introduced above. 

(i) (if) Suppose first that $j$ is nonpivotal. Let $p_m<j<p_{m+1}$. Since $j$
is inessential
\beq \label{eqv}
C_m[j] \in 
\cspan\{C_m[j+1], C_m[j+2], \ldots ,C_m[n]\}. \eeq
Let $(a_1,\ldots ,a_n)=\alpha_1R[1] + \ldots + \alpha_kR[k] \in X$.
Suppose that $a_i=0$ for $i=j+1,\ldots ,n$. Then, since $M$ is in rref,
we have $\alpha_{m+1}=\cdots = \alpha_k = 0$ and by (\ref{eqv}) above
we have $a_j=0$. Thus there is no vector in $X$ that has 
last nonzero component (from the left)
in column $j$. Thus $j\not\in \cR(X)$ and $s_j=H$.  

Now assume that $j$ is pivotal. Let $j=p_m$. Since $j$ is inessential
$$R_{n-j}[m] \in 
\cspan\{R_{n-j}[1], R_{n-j}[2], \ldots ,R_{n-j}[m-1]\}.$$
We have $R_{n-j}[m]=\alpha_1R_{n-j}[1]+\cdots +\alpha_{m-1}R_{n-j}[m-1]$. It
now follows that
$$R[m] - \left(\alpha_1R[1]+\cdots +\alpha_{m-1}R[m-1]\right)$$
has last nonzero component in column $j$. So $j\in \cR(X)$ and $s_j=H$.

The only if part is similar.

(ii) (if)
Let $j=p_m$. 
Let $(a_1,\ldots ,a_n)=\alpha_1R[1] + \ldots + \alpha_kR[k] \in X$.
Suppose that $a_i=0$ for $i=j+1,\ldots ,n$. Then, since $M$ is in rref,
we have $\alpha_{m+1}=\cdots = \alpha_k = 0$. Since $j$ is essential 
we have $\alpha_m = 0$ and so $a_j=0$. Thus no vector in $X$ has last nonzero component
in column $j$.
So $j\not\in \cR(X)$ and $s_j=U$.

The only if part is similar.

(iii) (if) Let $m$ satisfy $p_m<j<p_{m+1}$. Since $j$ is essential 
we can find a row vector   
$a=(a_1,\ldots ,a_m)$ such that $aC_m[j]\not=0$ and $aC_m[i]=0$ for 
$i=j+1,\ldots n$. It follows that $a_1R[1]+\cdots +a_mR[m]$ has last
nonzero component in column $j$ and so $j\in \cR(X)$. Thus $s_j=D$.

The only if part is similar. \eprf

{\bf{Proof of Theorem \ref{qbiane}}} Let $\Psi(X)=s_1s_2\cdots s_n$. We shall
prove by induction on $j$ that 
\beq \label{pr}
&s_1\cdots s_j \mbox{ is a Motzkin path from $(0,0)$
to $(j,t)$, where $t= \rank\;S_j$.}&\eeq 
Since $\rank\;S_n=0$ this will prove the result.
The claim is clearly true for $j=0$. Assume
the claim has been proved upto some $j\geq 0$ and we have built a Motzkin path
$P$ from $(0,0)$ to $(j,t)$, where $t=\rank\;S_j$.

Consider $C[j+1]$. The following cases arise:

(i) $j+1$ is pivotal and essential: By Theorem \ref{cqbiane}, 
the next point on $P$
is $(j+1,t+1)$. Clearly $\rank\;S_{j+1}=t+1$ (in both cases, $j$ pivotal and
$j$ nonpivotal). 

(ii) $j+1$ is pivotal and inessential: By Theorem \ref{cqbiane}, 
the next point on $P$
is $(j+1,t)$. Clearly $\rank\;S_{j+1}=t$. 
 
(iii) $j+1$ is nonpivotal and inessential: By Theorem \ref{cqbiane}, 
the next point on $P$
is $(j+1,t)$. Clearly $\rank\;S_{j+1}=t$. 
 
(iv) $j+1$ is nonpivotal and essential: Since $j+1$ is
essential, we have $\rank\;S_j \geq 1$ (in both cases, $j$ pivotal and $j$ nonpivotal).
By Theorem \ref{cqbiane}, 
the next point on $P$ is $(j+1,t-1)$ (note that $t-1\geq 0$). Clearly
$\rank\;S_{j+1}=t-1$.
That completes the inductive proof. \eprf
 
We now work towards the proof of Theorem \ref{main}. A matrix in rref is
said to be {\em{primary}} if all inessential columns are nonpivotal. 
A subspace $X$
is {\em{primary}} if $\row(X)$ is primary.  

\bt \label{bpb}
Let $P\in \cM(n)$. Then the number of primary subspaces $X\in B_q(n)$ 
with $\Psi(X)=P$ is given by $(q-1)^{|P|}w(P,q)$.
\et
\pf By induction on $n$, the case of $n=1$ being clear. Let $X\in B_q(n)$,
$M=\row(X)$, and write $\Psi(M)=s_1s_2\cdots s_n$. 
Let $t\geq 0$ be the largest
integer with $s_1=s_2=\cdots =s_t=U$. So $s_{t+1}=H \mbox{ or }D$.

From $s_1=\cdots =s_t=U$ we deduce that the first $t$ columns of $M$ 
are pivotal and essential and from (\ref{pr}) we have $\rank\;S_t = t$.

Now consider two cases (we continue to use the notation introduced above for the columns of $M$):

(i) $s_{t+1}=H$: So column $t+1$ of $M$ is inessential and by assumption is 
nonpivotal. Thus $S_{t+1}$ is $t\times (n-t-1)$ and  $\rank\;S_{t+1}$ is $t$ and $C_t[t+1]$ is in the
column space of $S_{t+1}$. The weight of the horizontal step
$s_{t+1}$ is $q^t$.

Let $P'$ be the Motzkin path $s_1\cdots s_ts_{t+2}\cdots s_n\in \cM(n-1)$.
We have $w(P,q)=q^tw(P',q)$. By induction hypothesis, 
the number of primary $(n-1)$-column rref's $M'$ with $\Psi(M')=P'$
is $(q-1)^{|P'|}w(P',q)$. 
Taking any such rref $M'$
and adding a new column $t+1$ with first $t$ coefficients arbitrary 
and others 0 (so there are $q^t$ choices for column $t+1$)
gives a $n$ column primary rref $N$ with $\Psi(N)=P$.
Since $|P'|=|P|$ the assertion follows. 

(ii) $s_{t+1}=D$: So column $t+1$ of $M$ is essential and nonpivotal. Thus $S_{t+1}$ is $t\times (n-t-1)$ and 
$\rank\;S_{t+1}=t-1$ and $C_t[t+1]$ is outside the column space of $S_{t+1}$.
The number of column vectors (with $t$ components) outside
the column space of $S_{t+1}$ is $q^t-q^{t-1}$. The weight of the down step $s_{t+1}$ is $q^{t-1} + q^t
+ \cdots + q^{2t-2}$.

Let $P'$ be the Motzkin path $s_1\cdots s_{t-1}s_{t+2}\cdots s_n \in \cM(n-2)$.
We have $w(P,q)=(q^{t-1}+q^t+\cdots +q^{2t-2})w(P',q)$.

Since $\rank\;S_{t+1}=t-1$, there is a unique nonzero row vector
$a=(a_1,\ldots ,a_t)$, with last nonzero coefficient (from the left) equal to 1
and with $aS_{t+1}=0$. Say the last nonzero coefficient of $a$ occurs in
column $j$, i.e., $a_{j+1}=\cdots = a_{t}=0$ and $a_j=1$. Then 
it is easy to see that removing row $j$
and columns $j,t+1$ of $M$ gives a $n-2$ column primary 
rref $M'$ with $\psi(M')=P'$.
Note that knowing $a$ and $M'$ we can recover row $j$ of $S_{t+1}$ from the
other rows and thus, we can recover $M$ except for column $t+1$. 

Observing that the $a$ in the paragraph above is unique 
we see by the induction hypothesis that the number of primary 
$M$ with $\Psi(M)=P$ is
$$(q^t-q^{t-1})(1+q+\cdots +q^{t-1})(q-1)^{|P'|}w(P',q)=(q-1)^{|P|}w(P,q),$$
since $|P|=|P'|+1$. That completes the proof. \eprf

\bex \label{pex}
Consider the following Motzkin path $P\in \cM(8)$. 
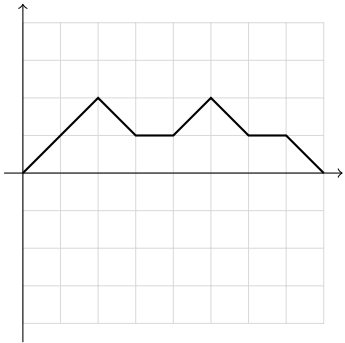
\begin{figure}[h]
    \centering
    \begin{tikzpicture}[scale=0.5]
        \draw[help lines, color=gray!30] (0,-4) grid (8,4);
        \draw[thick] (0,0) -- ++(1,1) -- ++(1,1) -- ++(1,-1) -- ++(1,0) -- ++(1,1) -- ++(1,-1) -- ++(1,0) -- ++(1,-1);
        \draw[->] (-0.5,0) -- (8.5,0) node[right] {};
        \draw[->] (0,-4.5) -- (0,4.5) node[above] {};
    \end{tikzpicture}
    \caption{$UUDHUDHD$}
    \label{fig:motzkin_path}
\end{figure}
An 8-column primary rref $M$ satisfies $\Psi(M)=P$ if and only if $\cL(M)= \{1,2,5\}$ and 
$$\rank\;S_1=1, \rank\;S_2=2, \rank\;S_3=1, \rank\;S_4=1, \rank\;S_5=2, \rank\;S_6=1, \rank\;S_7=1.$$
Such an example, which is used later, is given by the following matrix:
$$M=\left[ \ba{cccccccc}
1&0&a&0&0&0&0&0\\
0&1&b&c&0&d&e&e\\
0&0&0&0&1&0&f&f \ea \right],$$
where $a,d,e,f\in \F_q^*$ and $b,c\in \F_q$.

\eex
We have now given an interpretation to the term $(q-1)^{|P|}w(P,q)$ in identity (\ref{qbi}) above. We shall now give an interpretation to the  coefficient $\binom{n-2|P|}{k-|P|}$. Although the underlying idea is simple the details are somewhat involved.
We begin by recall a standard result from linear algebra and a crucial result from {\bf\cite{vv}}. 

Let $b,c\in \F_q^n$.
Set
$$\Gamma(b,c)=I + b^Tc,$$
where $I$ is the $n\times n$ identity matrix and $T$ denotes transpose. 
It is easily seen that
$\det(\Gamma(b,c))= 1+bc^T$ and so $\Gamma(b,c)$ is invertible if and 
only if $bc^T\not = -1$.
In this case the Sherman-Morrison-Woodbury formula gives
$$\Gamma(b,c)^{-1}= I - \frac{1}{1+bc^T}b^Tc.$$

The following result is from {\bf\cite{vv}}. For completeness
we include their proof.

\bl \label{vvl}
There is a bijection $\phi_n : \F_q^n \rar \F_q^n$ 
 such that $\Gamma(b,\phi_n(b))$ is nonsingular (i.e., 
$b\phi_n(b)^T\not= -1$) for all $b\in \F_q^n$.
\el
\pf We proceed by induction on $n$, beginning with the case $n=1$. For every $d\in \F_q,\;
d\not= 0$ define 
$$\mu_d:\F_q\rar \F_q$$
by $\mu_d(x)=1$ if $x=0$ and $\mu_d(x)=1 +dx^{-1}$ if $x\not=0$. Then $\mu_d$
is a bijection and $x\mu_d(x)\not= d$, for all $x\in \F_q$. Define 
$\phi_1 = \mu_{-1}$.

Now let $b=(b_1,\ldots ,b_{n+1})$. Write 
$\phi_n((b_1,\ldots ,b_n))=(c_1,\ldots ,c_n)$. Define 
$$\phi_{n+1}(b)=(c_1,\ldots ,c_n,\mu_{\alpha}(b_{n+1})),$$
where $\alpha= -1 - (b_1c_1+\cdots +b_nc_n)$. By induction hypothesis
$\alpha\not=0$, so $\mu_{\alpha}(b_{n+1})$ is defined and thus, since 
$x\mu_d(x) \not= d$ for all $x\in \F_q$, we have $b\phi(b)^T \not= -1$. \eprf

Note that, given $b\in \F_q^n$, both $\phi_n(b)$ and $\phi_n^{-1}(b)$
can be efficiently calculated recursively. This is important when we need to 
calculate the element covering a given subspace in the SCD
of $B_q(n)$.

Let $P\in \cM(n)$, $X\in B_q(n)$ with $\row(X)=M$ and $\Psi(M)=P$.
Let $j\in [n]$ be an inessential, pivotal column of $M$. We define a rref
$\dele(M,j)$ as follows. 
We shall use the notation for the rows, columns, and pivotal columns of $M$
introduced above.

Write $\cL(M)=\{p_1<p_2<\cdots <p_k\}$ and let $j=p_m$. 
We can picture $M$ as

\beq \label{delp}
\left[ \ba{ccc} A&0&N\\0&1&a\\0&0&B \ea\right],
\eeq

where $a$ is the row vector $R_{n-j}[m]$ and $\left[ \ba{c} N \\ a \ea \right]$ is the 
section at $j$. Note
that $a$ is in the row space of $N$, since $j$ is inessential.

Let $\rank\;N =s$ and let $R_{n-j}[i_1], R_{n-j}[i_2], \ldots ,R_{n-j}[i_s]$,
where $1\leq i_1 < i_2 <\cdots <i_s \leq m-1$, be the lexically first basis
of $N$, i.e., the first nonzero row of $N$ has index
$i_1$, the first row after row $i_1$ that is independent of row $i_1$
is row $i_2$ and so on. Let $N_L$ denote the submatrix of $N$ consisting of the
rows with indices $i_1,i_2,\ldots ,i_s$. 
There is a unique row vector $c=(c_1,c_2,\ldots ,c_s)$
such that
$$a=cN_L.$$
%
Let $(b_1,\ldots ,b_s)=\phi_s^{-1}((c_1,\cdots ,c_s))$.

Now define a column vector $d=(d_1,\ldots ,d_{m-1})^T$ as follows. First set
$d_{i_l}=b_l,\;l=1,\ldots ,s$. The other components of $d$ 
are defined by linearity.
Consider $d_u$, where $u\not\in\{i_1,\ldots ,i_s\}$.
The row $R_{n-j}[u]$ can be uniquely written as 
$$R_{n-j}[u]=\alpha_1R_{n-j}[i_1]+\cdots +\alpha_sR_{n-j}[i_s].$$
Set $d_u =\alpha_1b_1+\cdots +\alpha_sb_s$. It is clear that $d$ is in the column space
of $N$.

Now perform the following row operations on $M$: for $l=1,\ldots ,m-1$, multiply
row $m$ of $M$ by $d_l$ and add to row $l$. We get the following matrix

$$\left[ \ba{ccc} A&d&N'\\0&1&a\\0&0&B \ea\right],$$

Now delete row $m$ from the matrix above to get $\dele(M,j)$, pictured below.

\beq \label{delpp}
\left[ \ba{ccc} A&d&N'\\0&0&B \ea\right],
\eeq

Set $Y=\sspace(\dele(M,j))$.
\bl \label{delete}
Let $P, X, M, j, Y$ be as above. Then

(i) $Y\seq X$.

(ii) $\cL(Y) = \cL(X)\setminus \{j\}$.

(iii) The index sets of the independent rows of $N$ and $N'$ are the same.
In particular, the index set of the lexically first basis of $N'$ is the same as that of 
$N$.

(iv) $\Psi(Y)=P$.
\el
\pf Parts (i) and (ii) are clear.

(iii) Let $W$ be the row space of $N$.
Let $N_L'$ denote the submatrix of $N'$ 
consisting of the rows with indices $i_1,i_2,\ldots ,i_s$. Then we have
$$N_L'=\Gamma(\phi_s^{-1}(c), c)N_L.$$
Since $\Gamma(\phi_s^{-1}(c),c)$ is nonsingular, it follows that there
is a linear bijection $W\rar W$ such that, for $1\leq i \leq m-1$, row
$i$ of $N'$ is the image, under this map, of row $i$ of $N$.
The result follows.

(iv) Since column $j$ is inessential in $X$ and $j\in \cL(X)$ we have, by
Theorem \ref{cqbiane}, that $j\in \cR(X)$.
Clearly column $j$ is inessential in $Y$ (i.e., $d$ is in the column space of $N'$) and since $j\not\in \cL(Y)$ we have, by
Theorem \ref{cqbiane}, that $j\not\in \cR(Y)$. Thus $\cR(Y)
= \cR(X)\setminus \{j\}$. The result follows. \eprf

Let $P\in \cM(n)$, $X\in B_q(n)$ with $\row(X)=M$ and $\Psi(M)=P$.
Let $j\in [n]$ be an inessential, nonpivotal column of $M$. 
We shall now define a rref $\ins(M,j)$ as follows. We use
notation defined previously.

Write $\cL(M)=\{p_1<p_2<\cdots <p_k\}$ and let $p_m<j<p_{m+1}$. 
We can picture $M$ as

$$\left[ \ba{ccc} A&d&N\\0&0&B \ea\right],$$

where $N$ is the section at $j$ and $d$ is the column vector $C_m[j]$. Note
that $d$ is in the column space of $N$, since $j$ is inessential.

Let $\rank\;N=s$ and let $R_{n-j}[i_1], R_{n-j}[i_2], \ldots ,R_{n-j}[i_s]$,
where $1\leq i_1 < i_2 <\cdots <i_s \leq m$, be the lexically first basis
of $N$. Define the row vector  $b=(d_{i_1},d_{i_2},\ldots ,d_{i_s})$ and
let $\phi_s(b)=c=(c_1,c_2,\ldots ,c_s)$. Define a row vector
$$a=c\Gamma(b,\phi_s(b))^{-1}N_L.$$
Note that $a$ is in the row space of $N$.
Add to $M$ the following row vector with a pivotal 1 in column $j$ to get the matrix $M'$ pictured below

$$\left[ \ba{ccc} A&d&N\\0&1&a\\0&0&B \ea\right],$$

Now perform the following row operations on $M'$: for $l=1,\ldots ,m$, multiply
row $m+1$ of $M'$ by $d_l$ and subtract it from row $l$. We get the matrix $\ins(M,j)$
pictured below

$$\left[ \ba{ccc} A&0&N'\\0&1&a\\0&0&B \ea\right],$$

Set $Y=\sspace(\ins(M,j))$. We define $N_L$ and $N_L'$ as above.
\bl \label{insert}
Let $P, X, M, j, Y$ be as above. Then

(i) $X\seq Y$.

(ii) $\cL(Y) = \cL(X)\cup \{j\}$.

(iii) The index sets of the independent rows of $N$ and $N'$ are the same.
In particular, the index set of the lexically first basis of $N'$ is the same as that of 
$N$.

(iv) $\Psi(Y)=P$.
\el
\pf Parts (i) and (ii) are clear.

(iii) Write $N_L$ as
$$N_L=\Gamma(b,\phi_s(b))\Gamma(b,\phi_s(b))^{-1} N_L.$$
We can think of $\Gamma(b,\phi_s(b))\Gamma(b,\phi_s(b))^{-1} N_L$ as follows: 
start with the matrix
$\Gamma(b,\phi_s(b))^{-1} N_L$ and for $l=1,\ldots ,s$ add $d_{i_l}$ times $a$ to row
$l$ of $\Gamma(b,\phi_s(b))^{-1} N_L$. These opeartions will be undone by 
our operations on $M'$ and so we have
$$N_L'=\Gamma(b,\phi_s(b))^{-1} N_L.$$
Since $\Gamma(b,\phi(b))^{-1}$ is nonsingular, the result follows.

(iv) This is similar to the proof of part (iv) of Lemma \ref{delete}. 
\eprf

\bl \label{inverse}
Let $M$ be a rref with $n$ columns.

(i) Let $j\in [n]$ be an inessential, pivotal column of $M$. Then 
$$\ins(\dele(M,j),j)=M.$$

(ii) Let $j\in [n]$ be an inessential, nonpivotal column of $M$. Then 
$$\dele(\ins(M,j),j)=M.$$
\el
\pf
(i) We use the notation set up in defining $\dele(M,\;)$ and $\ins(M,\;)$.

Picture $M$ as in (\ref{delp}) and $\dele(M,j)$ as in (\ref{delpp}) and consider 
the vectors $a,b,c,d$ defined there. The index sets of the lexically first bases
in $N$ and $N'$ are the same and from Lemma \ref{delete} we have
$$N_L'=\Gamma(\phi_s^{-1}(c),c)N_L.$$

Now picture $\ins(\dele(M,j),j)$ as follows.

$$\left[ \ba{ccc} A&0&N''\\0&1&a'\\0&0&B \ea\right],$$

From Lemma \ref{insert} we have that the index sets of the lexically first bases of
$N'$ and $N''$ are the same and thus we have (note that $\Gamma$ is invoked with the same arguments $(\phi_s^{-1}(c),c)$ as above since the index sets of the lexically first bases of $N$ and $N'$ are the same)
\beqn
N_L'' &=& \Gamma(\phi_s^{-1}(c),c)^{-1}N_L'\\
      &=& \Gamma(\phi_s^{-1}(c),c)^{-1}\Gamma(\phi_s^{-1}(c),c)N_L\\
      &=& N_L.
\eeqn 
By linearity it follows that $N''=N$. Now
$$a'=c\Gamma(\phi_s^{-1}(c),c)^{-1}N_L' = cN_L=a,$$
completing the proof.

(ii) Similar to part (i). \eprf

We now define how to add/delete more than one inessential column.
Let $X\in B_q(n)$ with $\row(X)=M$. Define
\beqn
\set(X) &=& \{j\in [n]\;|\;\mbox{column $j$ of $M$ is inessential}\},\\
\subbset(X)&=& \set(X)\cap \cL(X).
\eeqn
Equivalently, $\set(X)=[n]\setminus (\cL(X)\Delta \cR(X))$ and 
$\subbset(X)=\cL(X)\cap \cR(X)$. 
Note that $X$ is primary if and only if $\subbset(X)=\emptyset$. Also
define $\set(M)=\set(X)$ and $\subbset(M)=\subbset(X)$.

Let $J\seq \subbset(X)$. List the elements of $J$ in increasing order as
$J=\{j_1 < j_2 <\cdots < j_t\}$. Define
$$\dele(M,J)=\dele(\cdots \dele(\dele(\dele(M,j_1),j_2),j_3)\cdots, j_t),$$
i.e., first delete $j_1$, then $j_2$, and so on up to $j_t$.   

Let $I\seq \set(X)\setminus \subbset(X)$. List the elements of $I$ in increasing order as
$J=\{j_1 < j_2 <\cdots < j_t\}$. Define
$$\ins(M,J)=\ins(\cdots \ins(\ins(\ins(M,j_t),j_{t-1}),j_{t-2})\cdots, j_1),$$
i.e., first insert $j_t$, then $j_{t-1}$, and so on up to $j_1$.   

\bl \label{Inverse}
Let $X\in B_q(n)$ with $\row(X)=M$. 

(i) Let $J\seq \subbset(X)$. Then $$\ins(\dele(M,J),J)=M.$$

(ii) Let $I\seq \set(X)\setminus \subbset(X)$. Then $$\dele(\ins(M,I),I)=M.$$
\el
\pf Follows from Lemma \ref{inverse}. \eprf

\bl \label{slb}
Let $P\in \cM(n)$ and $M$ be a rref with $\Psi(M)=P$. Let $j<i$ be
inessential columns of $M$ and let $j$ be nonpivotal. Let $N$ be the section
of $M$ at $j$, of size $m\times n-j$.

(i) Assume $i$ is nonpivotal and let $N'$ be the section at $j$ of $\ins(M,i)$.
Then $N'$ is also of size $m\times n-j$ and the index sets of the independent
rows of $N$ and $N'$ are the same. 

(ii) Assume $i$ is pivotal and let $N'$ be the section at $j$ of $\dele(M,i)$.
Then $N'$ is also of size $m\times n-j$ and the index sets of the independent
rows of $N$ and $N'$ are the same. 
\el
\pf (i) Clearly the sizes of $N,N'$ are the same. Picture $N,N'$ as follows

$$\left[ \ba{ccc} N_1&d&N_2 \ea\right]\;\;\;\left[ \ba{ccc} N_1&0&N_2'\ea\right],$$

where $d$ denotes column $i$ of $N$.
We have

(a) $d$ is in the column space of $N_2$.

(b) Let $W$ be the row space of the section at $i$ of $M$.
According to Lemma \ref{insert} there is a linear bijection $\Lambda:
W \rar W$ such that row $l$ of $N_2'$ is obtained
by applying $\Lambda$ to row $l$ of $N_2$.

The result follows from items (a), (b) above.

(ii) Similar to part (i). \eprf

The following result, along with Theorem \ref{bpb},  proves Theorem \ref{main}.
\bt \label{maind}
(i) Let $X\in B_q(n)$, with $\row(X)=M$ and $\Psi(X)=P$, be primary. Then
$$ \{\sspace(\ins(M,J))\;|\;J\seq \set(X)\}$$
is a symmetric Boolean subset of $B_q(n)$, with minimum rank $|P|$ and
maximum rank $n-|P|$.

(ii) We have the following SBD of $B_q(n)$
$$ B_q(n) = \coprod_P \coprod_M \;\{\sspace(\ins(M,J))\;|\;J\seq \set(M)\},$$
where $P$ varies over all Motzkin paths in $\cM(n)$ and $M$ varies over all
primary rref's with $\Psi(M)=P$.
\et
\pf
Write (i) $P=s_1s_2\cdots s_n$. Set $Z=\{i\in [n]\;|\; s_i=U\}$.

Let $I,J\seq \set(X)$ and put $M_1=\ins(M,I)$ and $M_2=\ins(M,J)$. We shall
show that $\sspace(M_1)\seq \sspace(M_2)$ if and only if $I\seq J$.
The only if part is clear since $\cL(M_1)=Z\cup I$ and $\cL(M_2)=Z\cup J$.

Now assume $I\seq J$. We shall prove by induction on $|J|$ that $\sspace(M_1)
\seq \sspace(M_2)$. If $|J|=0$ then there is nothing to prove. Assume we have
proved the result up to $|J|\leq k$. Let $|J|=k+1$. We may assume $|I|<|J|$.
Let $j=\min J$. If $j\not\in I$, then
$$\sspace(\ins(M,I))\seq \sspace(\ins(M,J\setminus \{j\}))\seq \sspace(\ins(M,J)),$$
where the first inclusion follows by the induction hypothesis and the second by
definition of $\ins(M,J)$ and Lemma \ref{insert}.

Now consider the case $j\in I$. Set 
$M_1'=\ins(M,I\setminus \{j\})$ and
$M_2'=\ins(M,J\setminus \{j\})$.
By the induction hypothesis,
$$\sspace(M_1')\seq \sspace(M_2').$$
We  have 
$$M_1 = \ins(M_1',j),\;\;\;M_2= \ins(M_2',j).$$
It follows
from our definition of insertion/deletion of inessential columns that the sections
at $j$ of $M_1'$ and $M_2'$ have the same size, say $m\times (n-j)$ and that the 
first $m$ components of column $j$ of $M_1'$ and $M_2'$ are identical. Call this column vector vector $d$.

Let the section at $j$ of $M_1',M_2'$ be denoted by $N_1',N_2'$ respectively.
Picture $M_1', M_2'$ as follows

$$\left[ \ba{ccc} A&d&N_1'\\0&0&B_1' \ea\right]\;\;\;
\left[ \ba{ccc} A&d&N_2'\\0&0&B_2' \ea\right],$$

and let $a_1, a_2$ be , respectively, the last rows of the sections at $j$ of
$M_1, M_2$. 

By Lemma \ref{slb} the index sets of the lexically first bases of $N_1'$ and $N_2'$
are the same (= the index set of the lexically first basis of the section at
$j$ of $M$). Since the first $m$ components in column $j$ of $M_1'$ and $M_2'$
are identical it follows from the definition of insertion that $a_1$ and $a_2$ are identical
linear combinations of the lexically first bases (having same index sets)
of $N_1'$ and $N_2'$. Since $\sspace(M_1')\seq \sspace(M_2')$ we see that
$a_2-a_1$ is in the row space of $B_2'$. Thus $\sspace(M_1)\seq \sspace(M_2)$.

(ii) For a rref $M$, $\dele(M,\subbset(M))$ is primary and, from Lemma \ref{Inverse}
$$M=\ins(\dele(M,\subbset(M)), \subbset(M)).$$
It remains to show that the union on the right hand side is disjoint.
Let $M_1,M_2$ be primary rref's and $J_1\seq \set(M_1), J_2 \seq \set(M_2)$.
Suppose $\sspace(\ins(M_1,J_1))=\sspace(\ins(M_2,J_2))$. From 
Lemma \ref{insert} $\Psi(M_1) = \Psi(M_2)$. Comparing pivotal columns
we get $J_1=J_2$. From Lemma \ref{Inverse} we now get
$$M_1=\dele(\ins(M_1,J_1), J_1)=\dele(\ins(M_2,J_2), J_2)=M_2,$$
completing the proof. \eprf

\bex \label{pexx}
Consider the primary rref $M$ from Example \ref{pex}. Then
$\set(M)=\{4,7\}$. We have
\newpage
\beqn
M &=&\left[ \ba{cccccccc}
1&0&a&0&0&0&0&0\\
0&1&b&c&0&d&e&e\\
0&0&0&0&1&0&f&f \ea \right],\\
\ins(M,7)&=&\left[ \ba{cccccccc}
1&0&a&0&0&0&0&0\\
0&1&b&c&0&d&0&\frac{e}{1+e\phi_1(e)}\\
0&0&0&0&1&0&0&\frac{f}{1+e\phi_1(e)}\\
0&0&0&0&0&0&1&\frac{e\phi_1(e)}{1+e\phi_1(e)}
\ea \right],\\
\ins(M,4)&=&\left[ \ba{cccccccc}
1&0&a&0&0&0&0&0\\
0&1&b&0&0&\frac{d}{1+c\phi_1(c)}&\frac{e}{1+c\phi_1(c)}&\frac{e}{1+c\phi_1(c)}\\
0&0&0&1&0&\frac{\phi_1(c)d}{1+c\phi_1(c)}&\frac{\phi_1(c)e}{1+c\phi_1(c)}&\frac{\phi_1(c)e}{1+c\phi_1(c)}\\
0&0&0&0&1&0&f&f
\ea \right],\\
\ins(M,\{4,7\})&=&\left[ \ba{cccccccc}
1&0&a&0&0&0&0&0\\
0&1&b&0&0&\frac{d}{1+c\phi_1(c)}&0&\frac{e}{(1+e\phi_1(e))(1+c\phi_1(c))}\\
0&0&0&1&0&\frac{\phi_1(c)d}{1+c\phi_1(c)}&0&\frac{\phi_1(c)e}{(1+e\phi_1(e))(1+c\phi_1(c))}\\
0&0&0&0&1&0&0&\frac{f}{1+e\phi_1(e)}\\
0&0&0&0&0&0&1&\frac{e\phi_1(e)}{1+e\phi_1(e)}
\ea \right].
\eeqn
\eex

As a corollary of the explicit SBD of $B_q(n)$ we get the following  algorithm for constructing an explicit SCD of $B_q(n)$, first given, with a different description, in {\bf\cite{vv}}.

For a finite set $J$ of positive integers, let $B(J)$ denote the poset of
all subsets of $J$. Clearly, we can transfer the explicit SCD of $B(|J|)$,
found by the bracketing procedure (see {\bf\cite{gk}}), to an explicit
SCD of $B(J)$ using the unique order preserving bijection
$\{1,2,\ldots ,|J|\} \rar J$. Given a symmetric chain $(x_1,\ldots ,x_h)$
we call $x_h$ the {\em top} element.

{\bf{Algorithm }} {\em (SCD of $B_q(n)$)}

{\bf Input:} $X\in B_q(n)$, given by a spanning set. 

{\bf Output:} $Y\in B_q(n)$ covering $X$ in an SCD of $B_q(n)$ or the statement
that $X$ is a top element. 

{\bf Method:}\begin{enumerate}

\item Form the matrix $N$ with rows given by the generating set of $X$.
   
\item Run the left to right Gauss elimination on $N$ to compute $\cL(X)$
and the (unique) rref $M$ with $\sspace(M)=X$.

\item Run the right to left Gauss elimination on $N$ to compute $\cR(X)$.

\item Set $J=[n]\setminus (\cL(X) \Delta \cR(X))$.

\item Set $I=\cL(X)\cap \cR(X)$.

\item IF $\left( I \mbox { is a top element in the SCD of } B(J) \right)$
THEN \\ Output ``$X$ is a top element". RETURN 

\item Let $I$ be covered by $I\cup \{j\}$ in the SCD of $B(J)$.\\
 Output $Y=\sspace(\ins(\dele(M,I), I\cup\{j\}))$. RETURN

\end{enumerate}

\noi The correctness of this algorithm follows directly from Theorem \ref{maind}. Applying this algorithm to Example \ref{pexx} above, we
see that $M$ is covered by $\ins(M,4)$, which in turn is covered by
$\ins(M,\{4,7\})$ in the SCD of $B_q(8)$ while $\ins(M,7)$ is a top element.

\noi {\bf{Remark}} All results of this section that do not 
explicitly refer to counting are valid over all fields. Let $F$ be a field
and let $B_F(n)$ denote the poset of subspaces of $F^n$ (row vectors) 
(this is an infinite graded poset of rank $n$ when $F$ is infinite). The definition
of the map $\Psi$, essential/inessential columns, 
and primary matrices are the same.
Lemma \ref{vvl} continues to hold. Insertion/deletion of inessential columns can be 
defined as before (with the same properties). The SBD of Theorem \ref{maind} 
continues to hold
$$ B_F(n) = \coprod_P \coprod_M \;\{\sspace(\ins(M,J))\;|\;J\seq \set(M)\},$$
where $P$ varies over all Motzkin paths in $\cM(n)$ and $M$ varies over all
primary rref's with $\Psi(M)=P$.

\section{A Problem}

We state a problem concerning a possible relation between the set theoretic
and linear analogs of symmetric chains in products of chains and $B_q(n)$.

For a finite set $S$, let $\C[S]$ denote the complex vector space with $S$ as
basis. 
Let $P$ be a finite graded poset of rank $n$ with rank function
$r: P\rar \{0,\ldots ,n \}$ and let
$P_i$ denote the set of elements of $P$ of rank $i$. 
We have a vector space direct sum
$$ \C[P]=\C[P_0]\oplus \C[P_1] \oplus \cdots \oplus \C[P_n].$$ 

An element $v\in \C[P]$ is {\em homogeneous} if $v\in \C[P_i]$ for some $i$. We
extend the definition of rank to nonzero homogeneous elements in the obvious way.
The {\em
up operator}  $\cU:\C[P]\rar \C[P]$ is defined, for $x\in P$, by
$\cU(x)= \sum_{y} y$,
where the sum is over all $y$ covering $x$. 
A {\em symmetric Jordan chain} in $\C[P]$ is a sequence
$$J=(v_1,\ldots ,v_h)$$
of nonzero homogeneous elements of $\C[P]$
such that $\cU(v_{i-1})=v_i$ for
$i < h$, $\cU(v_h)=0$, and $r(v_1) + r(v_h) = n$
if $h\geq 2$, or $2r(v_1)= n$ if $h=1$.
Note that the
elements of this sequence are linearly independent, being nonzero and of
different ranks. 
A {\em symmetric Jordan basis} (SJB) of $\C[P]$ is a basis of $\C[P]$
consisting of a disjoint union of symmetric Jordan chains
in $\C[P]$. We think of a SJB as a linear analog of a SCD.
Just like in the case of SCD's, we distinguish between existence and 
construction for SJB's. 

A graded poset may have a SCD but no SJB and, similarly, it may have
a SJB but no SCD. However, it it has both a SCD and a SJB we can think
of a possible relation between them. We now make a definition to discuss this.

Let $J$ be the symmetric Jordan chain in $\C[P]$ displayed above. We say
that a symmetric chain $C=(x_1,\ldots ,x_h)$ in $P$ is {\em supported} by
$J$ if, for all $i$,  $r(x_i)=r(v_i)$ and $x_i$ is in the support of $v_i$
(i.e., $x_i$ occurs with a nonzero coefficient when the vector $v_i$ is written
as a linear combination of the elements of $P$ of rank $r(x_i)$).

Let $J=\{J_1,\ldots ,J_t\}$, each $J_i$ a symmetric Jordan chain, be a SJB of 
$\C[P]$ and let $C=\{C_1,\ldots ,C_t\}$, each $C_i$ a symmetric chain, be a SCD
of $P$. We say that $J$ {\em supports} $C$ if there is a bijection
$\pi:C\rar J$ such that, for all $i$, $C_i$ is supported by $\pi(C_i)$. 

(i) Let $C(n)=\{0<1<\cdots <n\}$ denote the chain of length $n$. Consider
the chain product $C(n_1,\ldots,n_k)=C(n_1)\times \cdots \times C(n_k)$ (so
$B(n)$ is isomorphic to the chain product $C(1)\times \cdots \times C(1)$ ($n$
factors)). Proctor {\bf\cite{p}} used the $\mathfrak{sl(2,\C)}$ method to show the existence of a SJB of $C(n_1,\ldots ,n_k)$.
A simple visual algorithm to inductively construct a SCD for 
$C(n_1,\ldots ,n_k)$ was given in {\bf\cite{btk}}. This algorithm
was linearized in {\bf\cite{s2}} yielding an inductive construction of an explicit
SJB of $\C[C(n_1,\ldots ,n_k)]$. Since both these algorithms are based on the 
same underlying idea there is the question of the precise relation between them.
In particular, is it possible to develop a systematic theory that 
allows us to extract
the SCD constructed in {\bf\cite{btk}} 
from the SJB constructed in {\bf\cite{s2}} and thereby show that the SCD is
supported by the SJB. For instance, Example 3.4 in {\bf\cite{s2}} writes down
the SJB of $\C[B(4)]$ produced by the algorithm and it can be seen by inspection
that it supports the SCD of $B(4)$ given by the bracketing procedure.

This problem is in the same spirit as Rota's problem {\bf\cite{bb+}} (see last
paragraph of page 211)
of explaining the precise relation between the RSK bijection and its linear analog, the straightening
law. Leclerc and Thibon {\bf\cite{lt}} developed a theory to extract the RSK
bijection from its linear analog.

(ii) Now consider the case of the subspace lattice. Terwilliger {\bf\cite{t}} showed the existence of an orthogonal (under the standard inner product) SJB of $\C[B_q(n)]$ (together with a formula for the ratio
of the lengths of the successive vectors in this SJB). In {\bf\cite{s3}} a linear
algebraic interpretation of the Goldman-Rota identity (\ref{gr}) was given
and was then used to construct a canonical, orthogonal SJB. Does this SJB support a SCD. In particular,
does this SJB support the SCD constructed in {\bf\cite{vv}} and studied further
in the present paper.
Is it possible to develop a theory that allows us to extract
an SCD from the SJB constructed in {\bf\cite{s3}}. 

\begin{center}
{\bf{Acknowledgement}}
\end{center}
We thank Zachary Hamaker for helpful remarks on the paper.

\end{document}